\documentclass[reqno]{amsart}
\usepackage{amssymb}
\usepackage{bbm,amssymb,amsfonts,amsmath,amsopn,amsthm}
 \newtheorem{thm}{Theorem}[section]
 
 \newtheorem{lem}[thm]{Lemma}
 \newtheorem{prop}[thm]{Proposition}
 \theoremstyle{definition}
 
 \theoremstyle{rem}

 \numberwithin{equation}{section}

\begin{document}

%
%
%
%
%
%
%
%
%

\title[Global existence and decay estimates for viscoelastic kirchhoff equation ]
 {Global existence and decay estimates for the viscoelastic kirchhoff equation with a delay term}

\author{Noureddine Sebih}

\address{Department of Mathematics,	University Djilali Liabes of Sidi-Bel-Abbes, Sidi Bel Abbes 22000, Algeria}

\email{noureddine$_{-}$sebih$@y$ahoo.com }

\author{Abdelhamid Mohammed Djaouti}

\address{Department of Mathematics,	University Hassiba BenBouali of Chlef}

\email{djaouti$_{-}$abdelhamid@yahoo.fr;\\
	a.mohammeddjouti@univ-chlef.dz}

\author{Chafi Boudekhil}

\address{Department of computer science,	University of Sidi Bel Abbes 22000, Algeria}

\email{ Chafi$_{-}$b$@y$ahoo.com}



\date{January 1, 2004}


\begin{abstract}
In this paper, we consider a viscoelastic kirchhoff equation with a delay term in the internal feedback. By using the Faedo-Galarkin approximation method we   prove the well-posedness of the global solutions. Introducing suitable energy, we prove the general uniform decay results.
\end{abstract}

\maketitle
\section{Introduction}
In this paper we investigate the global existence and uniform decay rate of the energy for solutions to the nonlinear
viscoelastic kirchhof problem with delay term in the internal feedback.
{\small \begin{equation}\label{1.1}
\begin{aligned}
|u'(x,t)|^{\rho}u''(x,t)+\Delta^{2}u(x,t)-\Delta u''(x,t)-M(\|\nabla u\|^{2})\Delta u(x,t) \\
-\int_{0}^{t}h(t-s)\Delta^2u(x,s)ds+\mu_{1}g(u'(x,t))+\mu_{2}g(u'(x,t-\tau))=0, &~~x\in\Omega,~t>0,
\end{aligned}
\end{equation}}
\begin{equation}\label{1.2}
u(x,t)=\frac{\partial u(x,t)}{\partial n}=0,~~x\in\partial\Omega,~t>0,
\end{equation}
\begin{equation}\label{1.3}
u(x,0)=u_{0}(x),~~u'(x,0)=u_{1}(x),~~x\in\Omega,
\end{equation}
\begin{equation}\label{1.4}
u'(x,t-\tau)=f_0(x,t-\tau),~~x\in\Omega,~~~t\in(0,\tau),
\end{equation}
where $\Omega$ is a bounded domain in $\mathbb{R}^{n}$ with smooth boundary $\partial\Omega$. $\frac{\partial}{\partial n}$ represents the outward normal derivative on $\partial\Omega$. $\rho$, $\mu_{1}$, $\mu_{2}$ are three positive real numbers, $M\in\mathcal{C}^{1}(\mathbb{R}^{+})$, $h$ is a positive non increasing function defined on $\mathbb{R}^{+}$ which represents the kernel of the memory term and $g$ is an odd non-decreasing function of the class ${C}^{1}(\mathbb{R})$ which represents internal feedback. \\
In the absence of the delay term, many authors have investigated   problem (\ref{1.1}) and proved the stability, instability and the exponential decaying energy of the system  under   suitable assumptions, see for example \cite{Cav2,Cav22,7A,MS2,MS3,San}. In   paper \cite{Cav2}, the authors considered a related problem with strong damping
  \begin{equation*} 
|u'(x,t)|^{\rho}u''(x,t)-\Delta u(x,t)-\Delta u''(x,t)-\int_{0}^{t}h(t-s)\Delta u(x,s)ds-\gamma\Delta u_t=0.
\end{equation*} 
They obtained the global existence result for $\gamma>0$ and the uniform exponential decay of the energy for $\gamma>0$. Lately, the decay result has been extended by \cite{MS3} to the case $\gamma=0$.\\
In a recent work  \cite{MS2}, Messaoudi and Tatar studied the following problem:
 \begin{equation*} 
|u'(x,t)|^{\rho}u''(x,t)-\Delta u(x,t)-\Delta u''(x,t)-\int_{0}^{t}h(t-s)\Delta u(x,s)ds=b|u|^{p-2}u.
\end{equation*}
By introducing a new functional and using a potential well method, they obtained the global existence of solutions and the uniform decay of the energy where the initial data are   stable in a suitable set. Han and Wang  proved in \cite{7A} the global existence and the uniform decay for the  following nonlinear viscoelastic equation with damping:
\begin{equation*} 
|u'(x,t)|^{\rho}u''(x,t)-\Delta u(x,t)-\Delta u''(x,t)-\int_{0}^{t}h(t-s)\Delta u(x,s)ds+u'=0.
\end{equation*}
It is well known that delay effects often arise in many pratical problems because these phenomen depend not only on the present state but also on the past  history of the system. In recent years, the behavior of solutions for PDEs with time delay effects has become an active area of research, see \cite{26AAAA,FLI,NiPi,NiPi2,SHP,YZ} and the references therein. Datko proved in \cite{26AAAA} that a small delay in a boundary control is a source of instability. To stabilize a system involving delay terms, additional control terms will be necessary. In \cite{NiPi} Nicaise and Pignotti  considered the following wave equation with a linear damping and delay term inside the  domain
\begin{equation*} 
u_{tt}-\Delta u+\mu_1u_t+\mu_2u_t(t-\tau)=0.
\end{equation*}
The stability was proved in the case $0<\mu_1<\mu_{2}$. Kirane and Said Houari in \cite{KBS} investigated the following linear viscoelastic wave equation with a linear damping and delay term
\begin{equation*} 
u_{tt}-\Delta u+\int_0^tg(t-s)\Delta u(s)ds+\mu_1u_t+\mu_2u_t(t-\tau)=0.
\end{equation*}
They showed that its energy was exponentially decaying when $0<\mu_2<\mu_1$. For the plate equation with time delay term, Park  consider in \cite{SHP} the problem
\begin{equation*} 
u_{tt}+\Delta^2u-M(\|\nabla u\|^{2})\Delta u+\sigma(t)\int_0^tg(t-s)\Delta u(s)ds+a_0u_t+a_1u_t(t-\tau)=0,
\end{equation*}
which can be considered as an extensive weak viscoelastic plate equation with a linear time delay term. The author obtained a general decay result of energy by using suitable energy and Lyapunov functionals. Yang in  \cite{YZ}  studied initial boundary value problem of Euler-Bernoulli viscoelastic equation with a delay term in the internal feedbacks,
\begin{equation*} 
u_{tt}+\Delta^2u-\int_0^tg(t-s)\Delta^2u(s)ds+\mu_{1}u_t+\mu_{2}u_t(t-\tau)=0.
\end{equation*}
A global existence and uniform decay rates for the energy was proved. Recently \cite{FLI} showed the energy decay of solutions for the following nonlinear viscoelastic equation with a time delay term in the internal feedback
\begin{eqnarray*}
&&u_{tt}+\Delta^2u-div\left(F(\nabla u)\right)-\sigma(t)\int_0^tg(t-s)\Delta^2u(s)ds\\
&&\qquad\qquad+\mu_{1}|u_t|^{m-1}u_t+\mu_{2}|u_t(t-\tau)|^{m-1}u_t(t-\tau)=0.
\end{eqnarray*}
In the present paper, we devote our study to problem (\ref{1.1})-(\ref{1.4}). We will prove the global existence of weak solutions
and the uniform exponential decay of the energy for this problem by using Faedo-Galerkin method and the perturbed
energy method, respectively.
Our paper is organized as follows. In Section 2, we present the assumptions and main results. Section 3 we prove our main results. 

\section{Assumptions and main result}

Let us consider the Hilbert space $L^{2}(\Omega)$ endowed with the inner product $(,)$ and the corresponding norm $\|\cdot\|$.
We also consider the sobolev space $H^{2}_{0}(\Omega)$ endowed with the scalar product
\begin{equation*}
(u,v)_{H^{2}_{0}(\Omega)}=(\Delta u, \Delta v).
\end{equation*}
We define for all $1\leq p<\infty$ and $u\in L^p(\Omega)$,
\begin{equation*}
\|u\|^{p}_{p}=\int_{\Omega}|u(x)|^{p}dx,~~\text{and}~~\|u\|=\|u\|_{2}.
\end{equation*}
We introduce as in \cite{NiPi} a new variable
\begin{equation*}
z(x,\varrho,t)=u_{t}(x,t-\tau\varrho),~~x\in\Omega,~\varrho\in(0,1),~t>0.
\end{equation*}
Then we have
\begin{equation*}
\tau z_{t}(x,\varrho,t)+z_{\varrho}(x,\varrho,t)=0.
\end{equation*}
Therefore problem (\ref{1.1})-(\ref{1.4}) is equivalent to
\begin{equation}\label{2.01}
\begin{aligned}
|u'(x,t)|^{\rho}u''(x,t)+\Delta^{2}u(x,t)-\Delta u''(x,t)-M(\|\nabla u\|^{2})\Delta u(x,t) \\
-\int_{0}^{t}h(t-s)\Delta^2u(x,s)ds+\mu_{1}g(u'(x,t))+\mu_{2}g(z(x,1,t))=0,
\end{aligned}
\end{equation}
\begin{equation}\label{2.02}
\tau z_{t}(x,\varrho,t)+z_{\varrho}(x,\varrho,t)=0,~~x\in\Omega,~\varrho\in(0,1),~t>0,
\end{equation}
\begin{equation}\label{2.03}
z(x,0,t)=u'(x,t),~~x\in\Omega,~t>0,
\end{equation}
\begin{equation}\label{2.04}
u(x,t)=\frac{\partial u(x,t)}{\partial n}=0,~~x\in\partial\Omega,~t>0,
\end{equation}
\begin{equation}\label{2.05}
u(x,0)=u_{0}(x),~~u'(x,0)=u_{1}(x),~~x\in\Omega,
\end{equation}
\begin{equation}\label{2.06}
z(x,\varrho,0)=f_{0}(x,-\varrho\tau),~~x\in\Omega,~\varrho\in(0,1).
\end{equation}
To state and prove our result, we assume  the following hypothesis
\begin{itemize}
	\item[(A1)] Assume that $\rho$ satisfies
	\begin{gather*}
	0<\rho\leq \frac{2}{n-2}\quad \text{if } n\geq 3,~~~~0<\rho< \infty\quad  \text{if } n=1,2;
	\end{gather*}
	\item[(A2)] Assume that  $M\in C^1(\mathbb{R}_+)$ satisfies
	\begin{equation*} 
	\begin{gathered}
	\exists m_0>0,~~M(\lambda)\geq m_0,~~\forall~\lambda\geq0.\\
	\exists \gamma,\delta,~~M(\lambda)\leq\delta\lambda^\gamma,~~\forall~\lambda\geq0.\\
	\exists\alpha,\beta,~~|M'(\lambda)|\leq\beta\lambda^\alpha,~~\forall~~\lambda\geq0.
	\end{gathered}
	\end{equation*}
	\item[(A3)] The kernel function $h:{\mathbb{R}_+}\to {\mathbb{R}_+^*}$ is a bounded
	$C^1$ function such that
	\begin{equation*} 
	1-\int_0^{\infty}h(s)\,ds=\beta_1>0,
	\end{equation*}
	and we assume that there exist a positive constant $\zeta $ satisfying
	\begin{equation*} 
	h'(t)\leq-\zeta h(t),~~t\geq0.
	\end{equation*}
	\item[(A4)] $g:\mathbb{R}\to \mathbb{R}$ is an odd non decreasing function of class $C^1$ such that there exist $c_{1}$,$\alpha_{1}$, $\alpha_{2}$ positives
	satisfying
	\begin{equation*} 
	\begin{gathered}
	|g'(s)|\leq c_1,~~\forall~s\in\mathbb{R},\\
	\alpha_1sg(s)\leq G(s)\leq\alpha_2sg(s),~~\forall~s\in\mathbb{R},
	\end{gathered}
	\end{equation*}
	where $G(s)=\int_0^sg(t)dt$, $\lim_{s\rightarrow+\infty}g(s)=+\infty$ and $\alpha_2\mu_2\leq\alpha_1\mu_1$.
\end{itemize}
First We state some lemmas which will be used in the next sections.
\begin{lem}\label{2.11}  For $h\in C^1([0,+\infty[,\mathbb{R})$ and $\varphi\in C^1(0,T;L^2(\Omega))$, we have
	\begin{equation*}
	\begin{array}{lll}
	&&\int_{\Omega}\int_{0}^{t}h(t-s)\varphi(x,s)\varphi'(x,t)dsdx\\
	&&=-\frac{1}{2}h(t)\|\varphi(t)\|^{2}+\frac{1}{2}(h'\square\varphi)(t)+\frac{1d}{2dt}\left[\left(\int_{0}^{t}h(s)ds\right)\|\varphi(t)\|^2-(h\square\varphi)(t)\right],
	\end{array}	
	\end{equation*}
	where $(h\square\varphi)(t)=\int_{0}^{t}h(t-s)\|\varphi(s)-\varphi(t)\|^{2}ds$.
\end{lem}

\begin{lem}\label{2.12}
	Let $\Phi$ is a convex function of class $C^1(\mathbb{R})$. The Legendre transformation of $\Phi$ is defined as follows $$\Phi^*(s)=\sup_{t\in\mathbb{R}}(st-\Phi(t)).$$
	If $\Phi'$ is an odd and $\lim_{s\rightarrow+\infty}\Phi'(s)=+\infty$, then
	\begin{equation*}
	\Phi^*(s)=s(\Phi')^{-1}(s)-\Phi\left((\Phi')^{-1}(s)\right),~~\forall~s\in\mathbb{R},
	\end{equation*}
	and satisfies the inequality
	\begin{equation*}
	st\leq\Phi^{*}(s)+\Phi(t),~~\forall~s,t\in\mathbb{R}.
	\end{equation*}
\end{lem}

The energy associated with problem (\ref{2.01})-(\ref{2.06}) is given by
\begin{equation}\label{2.13}
\begin{split}
E(t)&=\frac{1}{\rho+2}\|u'(t)\|^{\rho+2}_{\rho+2}+\frac{1}{2}\|\Delta u(t)\|_2^2+\frac{1}{2}\|\nabla u'(t)\|^{2}\\
&\quad+\frac{1}{2}\widehat{M}(\|\nabla u(t)\|^{2}) -\frac{1}{2}\left(\int_0^t h(s)\,ds\right)\|\Delta u(t)\|^{2}\\
&\quad +\frac{1}{2}(h\square\Delta u)(t)+\xi\int_\Omega\int_0^1G(z(x,\varrho,t))\,d\varrho\,dx,
\end{split}
\end{equation}
where $\widehat{M}(\lambda)=\int_{0}^{\lambda}M(t)dt$ and $\xi$ is a positive constant such that
\begin{equation*}
\tau \frac{\mu_2(1-\alpha_1)}{\alpha_1}<\xi<\tau\frac{\mu_1-\alpha_2 \mu_2}{\alpha_2}.
\end{equation*}
\begin{thm}\label{2.14}
	Let $u_0 \in H_0^2(\Omega)\cap H^3(\Omega)$, $u_1 \in H_0^2(\Omega)$ and
	$ f_0\in H_0^1(\Omega,H^1(0,1))$ satisfy the compatibility condition
	$f(\cdot,0)=u_1$.
	Assume that {\rm (A1)-(A4)} hold. Then (\ref{2.01})-(\ref{2.06}) admits a weak solution
{\small 	\begin{equation*}
	u\in L^{\infty}([0,\infty);H_0^2(\Omega)\cap H^{3}(\Omega)),~~u'\in L^{\infty}([0,\infty); H_0^2(\Omega)),~~u''\in L^2([0,\infty);\ H_0^1(\Omega)).
	\end{equation*}}
{\small 	\begin{equation*}
	z\in L^\infty([0,\infty); H_0^1(\Omega\times(0,1))),~~z'\in L^\infty([0,\infty); L^2(\Omega\times(0,1))),
	\end{equation*}}
{\small 	\begin{equation*}
	~~G(z(x,\varrho,t))\in L^{\infty}([0,\infty);L^1(\Omega\times(0,1))).
	\end{equation*}}
	Moreover, if $E(0)$ is positive and bounded, then for every $t_0>0$, there exist positive constants k and K such that the energy defined
	by (2.11) possesses the following decay:
	\begin{equation}\label{2.15}
	E(t)\leq Ke^{-kt},~~\forall t\geq t_0.
	\end{equation}
\end{thm}

\section{Proof of the main result}

We will divide the proof into two steps: in the first step, we will use the Faedo-Galerkin method to prove the existence of global
solutions, where the second  step  is devoted to proving the uniform decay of the energy by the perturbed energy method.\\\\
\textbf{Step 1}. Existence of weak solutions. \\
Let $T > 0$ be fixed and let $(w_i)_{i\in\mathbb{N}^{*}}$ be an orthogonal basis of $H^3(\Omega)\cap H^{2}_{0}(\Omega)$ with $w_i$ being the eigenfunctions of the bi-Laplacien operator subject to the boundary condition
\begin{equation*}
\Delta^2 w_i= \lambda_i w_i,\quad \text{in } \Omega,\quad w_i=\frac{\partial w_i}{\partial n} =0\quad \text{in } \partial\Omega.
\end{equation*}
By the linear elliptic operator theory  described in \cite{ZHeng}, we have $w_j\in H^m(\Omega)\cap H^{2}_{0}(\Omega)$, $m\in\mathbb{N}$. Now we denote by $W_{k}=\text{span}\{w_1,w_2,...,w_k\}$ the subspace generated by the first $k$ vectors of the basis $(w_i)_{i\in\mathbb{N}^{*}}$. By normalization, we
get $\|w_i\|=1$. Now we definie for $1\leq i\leq k$ the sequence as follows $\varphi_i(x,0)=w_i(x)$. Then we may extend $\varphi_i(x,0)$ by $\varphi_i(x,\varrho)$ over $L^2(\Omega\times]0,1[)$ and denote $Z_k$ the space generated by ${\varphi_1,...,\varphi_k}$. For any given integer $k$, we consider the approximate solution $(u_k,z_k)$
\begin{equation*}
u_k(x,t)=\sum_{i=1}^k g_{ik }(t)w_i(x), \quad
z_k(x,\varrho,t)=\sum_{i=1}^k h_{ik}(t)\varphi_i(x,\varrho),
\end{equation*}
which satisfies
{\small \begin{gather}\label{3.001}
\begin{gathered}
\begin{aligned}
&(|u'_k(t)|^{\rho}u''_k(t),w_i) +(\Delta u_k(t),\Delta w_i)
+(\nabla u''_k(t),\nabla w_i)+M(\|\nabla u_k(t)\|^{2})(\nabla u_k(t),\nabla w_i)\\
&-\int_0^t h(t-s)(\Delta u_k(s),\Delta w_{i})ds+ \mu_1(g(u'_k(t)),w_i)+\mu_2(g(z_k(.,1,t)),w_i)=0,
\end{aligned}\\
z_k(x,0,t)=u'_k(x,t),
\end{gathered}
\end{gather}}
\begin{equation}\label{3.002}
( \tau z_{kt}(t)+ z_{k\varrho}(t),\varphi_i)_{L^2(\Omega\times]0,1[)}=0 
\end{equation}
\begin{equation}\label{3.003}
u_k(0)=u_{0k},~~u'_k(0)=u_{1k},~~z_{k}(0)=z_{0k}
\end{equation}
where $i=\overline{1,k}$
\begin{equation*}
u_{0k}=\sum_{i=1}^k(u_0,w_i)w_i,~~u_{1k}=\sum_{i=1}^k(u_1,w_i)w_i,~~z_{0k}=\sum_{i=1}^k(f_0,\varphi_i)_{L^2(\Omega\times]0,1[)}\varphi_i
\end{equation*}
and for $k\to +\infty$
\begin{equation}\label{441}
\begin{array}{lll}
u_{0k}\to u_0~~&\text{in }&H_0^2(\Omega)\cap H^3(\Omega),\\
u_{1k}\to u_1~~&\text{in }&H_0^2(\Omega),\\
z_{0k}\to f_0~~&\text{in }&H_0^1 (\Omega,H^1(0,1)).
\end{array}
\end{equation}
Taking account of assumption (A1), $H^1_0(\Omega)\hookrightarrow L^{2(\rho+1)}(\Omega)$, then $u''_{k}\in L^{2(\rho+1)}(\Omega)$, $|u'_k|^\rho\in L^{\frac{2(\rho+1)}{\rho}}(\Omega)$ and $w_i\in L^2(\Omega)$, from the generalized H\"{o}lder inequality, the nonlinear term in (3.1)
$$(|u'_k(t)|^{\rho}u''_k(t),w_i)=\int_{\Omega}|u'_k|^\rho u''_k w_i dx\leq\|u'_k\|^\rho_{2(\rho+1)}\|u''_k\|_{2(\rho+1)}\|w_i\| $$
make sens. According to the standard of ordinary differential equations theory, the finite dimensional problem (\ref{3.001})-(\ref{3.003}) has a solution $(g_{ik},h_{ik})$ defined on $[0,t_{k}[$. Then we can obtain an approximate solution $u_k$ and $z_k$ of (\ref{3.001})-(\ref{3.003}) in $W_k$ and $Z_k$ respectively over $[0,t_k[$. Moreover, the solution can be extended to $[0,T]$ for any given $T$ by the first estimates below.\\
Now we derive the first estimate. Multiplying (\ref{3.001}) by $g'_{ik}(t)$ and summing with respect to $i$, we conclude that
\begin{eqnarray}
\nonumber&&\frac{d}{dt}\left[\frac{1}{\rho+2}\|u'_{k}(t)\|^{\rho+2}_{\rho+2}+\frac{1}{2}\|\Delta u_k(t)\|^{2}+\frac{1}{2}\|\nabla u'_{k}(t)\|^{2}+\frac{1}{2}\widehat{M}(\|\nabla u_k(t)\|^{2})\right]\label{0001}\\
&&\quad-\int_{0}^t h(t-s)(\Delta u_{k}(s),\Delta u'_{k}(t))ds\\
\nonumber&&\quad+\mu_1\int_\Omega u'_k(x,t)g(u'_k(x,t))dx+\mu_2\int_\Omega u'_k(x,t)g(z_k(x,1,t))dx=0.
\end{eqnarray}
Applying lemma 2.1 with $\varphi=\Delta u_k$,   (\ref{0001}) become
\begin{equation}\label{0002}
\begin{split}
&\frac{d}{dt}\Big[\frac{1}{\rho+2}\|u'_{k}(t)\|^{\rho+2}_{\rho+2}+\frac{1}{2}\|\Delta u_k(t)\|^{2}+\frac{1}{2}\|\nabla u'_{k}(t)\|^{2}+\frac{1}{2}\widehat{M}(\|\nabla u_k(t)\|^{2})\\
&\quad-\frac{1}{2}\left(\int_0^t h(s)\,ds\right)\|\Delta u_k(t)\|^{2}+\frac{1}{2}(h\square\Delta u_k)(t)\Big] \\
&=-\mu_1\int_\Omega u'_k(x,t)g(u'_k(x,t))dx-\mu_2\int_\Omega u'_k(x,t)g(z_k(x,1,t))dx\\
&\quad~~-\frac{1}{2}h(t)\|\Delta u_k(t)\|^2+\frac{1}{2}(h'\square\Delta u_{k})(t).
\end{split}
\end{equation}
We  multiply equation (\ref{3.002}) by $\xi g(z(x,\varrho,t))$ and integrating over $\Omega\times(0,1)$, we obtain
\begin{align*}
\xi \int_\Omega\int_0^1 z_{kt}(x,\varrho,t)g(z_k(x,\varrho,t))d\varrho dx
&= -\frac{\xi}{\tau} \int_\Omega\int_0^1 z_{k\varrho}(x,\varrho,t)g(z_k(x,\varrho,t))d\rho dx\\
&= -\frac{\xi}{\tau}\int_\Omega\int_0^1\frac{\partial}{\partial\varrho}\Big(G(z_k(x,\varrho,t))\Big)d\varrho dx .
\end{align*}
Hence
{\small \begin{equation}\label{9}
\xi \frac{d}{dt}\int_\Omega\int_0^1 G(z_k(x,\varrho,t))d\varrho dx =-\frac{\xi}{\tau} \int_\Omega G(z_k(x,1,t))dx+\frac{\xi}{\tau} \int_\Omega G(u'_k(x,t))dx.
\end{equation}}
Combining (\ref{0002}) and (\ref{9}), we obtain
\begin{equation}\label{0003}
\begin{split}
E'_k(t)
&=-\mu_1\int_\Omega u'_k(x,t)g(u'_k(x,t))dx-\mu_2\int_\Omega u'_k(x,t)g(z_k(x,1,t))dx\\
&  \quad\quad-\frac{1}{2}h(t)\|\Delta u_k(t)\|^2+\frac{1}{2}(h'\square\Delta u_{k})(t)-\frac{\xi}{\tau} \int_\Omega G(z_k(x,1,t))dx\\
&  \quad\quad+\frac{\xi}{\tau} \int_\Omega G(u'_k(x,t))dx.
\end{split}
\end{equation}
From assumption (A4),we knew that $G$ is a convex function of classe $C^{2}$, $G'=g$ is an odd and $\lim_{s\rightarrow+\infty}G'(s)=+\infty$, then by lemma 2.2, we deduce
$$G^*(s)=sg^{-1}(s)-G(g^{-1}(s)),~~\forall s\in\mathbb{R}.$$
Applying these equality with  $s=g(z_k(x,1,t))$, we obtain
\begin{equation*}
G^*(g(z_k(x,1,t)))=z_k(x,1,t)g(z_k(x,1,t))-G(z_k(x,1,t)),
\end{equation*}
By using inequality in lemma 2.2 together with (A4) for $s=g(z_k(x,1,t))$, $t=-u'_k(x,t)$ and $G$ is even function, we get
\begin{equation}\label{0002,}
\begin{split}
&-u'_k(x,t)g(z_k(x,1,t))\\
&\quad\quad\leq G^*(g(z_k(x,1,t)))+G(-u'_k(x,t)),\\
&\quad\quad= z_k(x,1,t)g(z_k(x,1,t))-G(z_k(x,1,t))+G(u'_k(x,t)),\\
&\quad\quad\leq (1-\alpha_{1})z_{k}(x,1,t)g(z_k(x,1,t))+\alpha_2u'_kg(u'_k).
\end{split}
\end{equation}
From (\ref{0003}), (\ref{0002,}) and assumption (A4), we have
\begin{equation}\label{0008}
\begin{split}
E'_k(t)
&\leq-\frac{1}{2}h(t)\|\Delta u_k(t)\|^2+\frac{1}{2}(h'\square\Delta u_{k})(t)\\
&\quad-\Big(\mu_1-\frac{\xi\alpha_2}{\tau}-\mu_2\alpha_2\Big)\int_\Omega u'_k(x,t)g(u'_k(x,t))dx\\
&\quad-\Big(\mu_1-\frac{\xi\alpha_2}{\tau}-\mu_2\alpha_2\Big)\int_\Omega z_k(x,1,t)g(z_k(x,1,t))dx.
\end{split}
\end{equation}
Integrating (\ref{0008}) over $(0,t)$ and using assumption $(A3)$, we conclude that
\begin{equation}\label{0005}
E_k(t)+\theta_{1}\int_{0}^{t}\int_\Omega u'_k(x,t)g(u'_k(x,t))dx+\theta_{2}\int_{0}^{t}\int_\Omega z_k(x,1,t)g(z_k(x,1,t))dx\leq C_{1}.
\end{equation}
where
\begin{equation*}
\theta_1=\Big(\mu_1-\frac{\xi\alpha_2}{\tau}-\mu_2\alpha_2\Big),~~\theta_2=\Big(\mu_1-\frac{\xi\alpha_2}{\tau}-\mu_2\alpha_2\Big),
\end{equation*}
and $C_{1}$ is a positive constant depending only on $\|u_{0}\|_{H_0^2(\Omega)}$, $\|u_1\|_{H_0^1(\Omega)}$ and $\|f_0\|_{L^2(\Omega\times(0,1))}$.
Noting (A2), (A3) and (\ref{0005}), we obtain the first estimate
\begin{eqnarray}
\nonumber&&\|u'_k(t)\|^{\rho+2}_{\rho+2}+\|\Delta u_k(t)\|^2+\|\nabla u'_k(t)\|^{2}\\
&&\label{1023}+\|\nabla u_k(t)\|^{2}+(h\square\Delta u_k)(t)+\xi\int_\Omega\int_0^1G(z_k(x,\varrho,t))\,d\varrho\,dx \\
\nonumber&&+\int_{0}^{t}\int_\Omega u'_k(x,s)g(u'_k(x,s))dxds+\int_{0}^{t}\int_\Omega z_k(x,1,s)g(z_k(x,1,s))dxds\leq C_{2},
\end{eqnarray}
where $C_{2}$ is a positive constant depending only on $\|u_{0}\|_{H_0^2(\Omega)}$, $\|u_1\|_{H_0^1(\Omega)}$, $\|f_0\|_{L^2(\Omega\times(0,1))}$, $m_0$, $\beta_{1}$, $\xi$, $\tau$, $\theta_{1}$ and $\theta_{2}$.\\
It follows from (\ref{1023}) that
\begin{eqnarray*}
&&u_k \text{ is uniformly bounded in }L^{\infty}(0,T;H_0^2(\Omega)),\\
&&u'_k \text{ is uniformly bounded in }L^{\infty}(0,T;H_0^1(\Omega)),\\
&&G(z_k(x,\rho,t)) \text{ is uniformly bounded in }L^{\infty}(0,T;L^1(\Omega\times(0,1))), \\
&&u'_kg(u'_k) \text{ is uniformly bounded in }L^1(\Omega\times(0,T)), \\
&&z_k(.,1,.)g(z_k(.,1,.)) \text{ is uniformly bounded in }L^1(\Omega\times(0,T))).
\end{eqnarray*}
Then, we derive the second estimate. Substituting $w_{i}$ by $-\Delta w_{i}$ in (\ref{3.001}), multiplying by $g'_{ik}$ and then summing with respect to $i$, it holds that
\begin{equation*}
\begin{split}
&\frac{1}{2}\frac{d}{dt}\left[\|\nabla \Delta u_k(t)\|^2+\|\Delta u'_k(t)\|^2\right]
-\int_\Omega|u'_k(x,t)|^\rho u''_k(x,t)\Delta u'_k(x,t) dx\\
&+\mu_1\int_\Omega |\nabla u'_k(x,t)|^2g'(u_k(x,t))dx+\mu_2\int_\Omega g'(z_k(x,1,t))\nabla u'_k(x,t)\nabla z_k(x,1,t)dx \\
&\qquad-\int_0^th(t-s) \int_{\Omega}\nabla\Delta u_k(x,s)\nabla\Delta u'_k(x,t)dxds\\
&\qquad+\int_{\Omega}M(\|\nabla u_k(t)\|^{2})\Delta u_k(x,t)\Delta u'_k(x,t)dx=0
\end{split}
\end{equation*}
Using lemma 2.1 with $\varphi=\nabla\Delta u_k$, we have
{\small \begin{equation} \label{20}
\begin{split}
&\frac{1}{2}\frac{d}{dt}\left[\left(1-\int_0^th(s)ds\right)\|\nabla \Delta u_k\|^2
+\|\Delta u'_k\|^2+M(\|\nabla u_k\|^2)\|\Delta u_{k}\|^{2}+h\square\nabla\Delta u_{k}(t)\right]\\
&\quad-\int_\Omega|u'_k|^\rho u''_k\Delta u'_k dx-M'(\|\nabla u_k\|^{2})\|\Delta u_{k}\|^{2}(\nabla u_k,\nabla u'_k)\\
&\quad+\mu_1\int_\Omega |\nabla u'_k |^2g'(u_k)dx+\mu_2\int_\Omega g'(z_k(x,1,t))\nabla u'_k\nabla z_k(x,1,t) dx \\
&=-\frac{1}{2}h(t)\|\nabla\Delta u_{k}\|^{2}+\frac{1}{2}(h'\square\nabla\Delta u_k)(t).
\end{split}
\end{equation}}
Substituting $\varphi_i$ by $-\Delta\varphi_i$ in (\ref{3.002}), multiplying by $h_{ik}$, summing with respect to $i$ and integratin over $\varrho\in(0,1)$, it follows that
\begin{equation} \label{21}
\frac{\tau}{2}\frac{d}{dt}\|\nabla z_k(t)\|^2_{L^2(\Omega\times(0,1))}+\frac{1}{2}\|\nabla z_k(.,1,t)\|^2-\frac{1}{2}\|\nabla u'_k(t)\|^2=0.
\end{equation}
Using the Green's formula, we get
\begin{eqnarray}\label{22}
 -\int_\Omega|u'_k|^{\rho}u''_k\Delta u'_kdx&=&\frac{d}{dt}\int_\Omega|u'_k(x,t)|^{\rho}|\nabla u'_k(x,t)|^2dx\\
&&\quad-\int_\Omega|u'_k(x,t)|^{\rho}\nabla u''_k(x,t)\nabla u'_k(x,t)dx.
\end{eqnarray}
Combining (\ref{20}) to (\ref{22}), we obtain
{\small \begin{eqnarray}
\nonumber&&\frac{1}{2}\frac{d}{dt}\Big[\left(1-\int_0^th(s)ds\right)\|\nabla \Delta u_k\|^2
+\|\Delta u'_k\|^2+M(\|\nabla u_k\|^2)\|\Delta u_{k}\|^{2}\\
\nonumber&&\quad+h\square\nabla\Delta u_{k}(t)+2\int_\Omega|u'_k|^{\rho}|\nabla u'_k|^2dx+\tau\|\nabla z_k(t)\|^2_{L^2(\Omega\times(0,1))}\Big]\\
&&\label{w15}\quad+\frac{1}{2}\|\nabla z_k(.,1,t)\|^2+\mu_1\int_\Omega |\nabla u'_k |^2g'(u_k)dx\\
\nonumber&&=\int_\Omega|u'_k|^{\rho}\nabla u''_k\nabla u'_kdx-\mu_2\int_\Omega g'(z_k(x,1,t))\nabla u'_k\nabla z_k(x,1,t)dx
+\frac{1}{2}\|\nabla u'_k(t)\|^2\\
\nonumber&&\quad+M'(\|\nabla u_k\|^{2})\|\Delta u_{k}\|^{2}(\nabla u_k,\nabla u'_k)-\frac{1}{2}h(t)\|\nabla\Delta u_{k}\|^{2}+\frac{1}{2}(h'\square\nabla\Delta u_k)(t)
\end{eqnarray}}
From Young inequality we have, for all $\eta>0$, that
\begin{equation*}
ab\leq\eta a^2+\frac{b^2}{4\eta}, ~~\text{where}~~a,b\in\mathbb{R}^{*}_{+}.
\end{equation*}
Assumption $(A2)$, (\ref{1023}) and using Young's inequality with $\eta=1/2$, we obtain
\begin{eqnarray*}
&&M'(\|\nabla u_k\|^{2})\|\Delta u_{k}\|^{2}(\nabla u_k,\nabla u'_k)\\
&&\qquad\leq M'(\|\nabla u_k\|^{2})\|\Delta u_{k}\|^{2}\|\nabla u_k\|\|\nabla u'_k\|\notag \\
&&\qquad\leq\frac{1}{2}\left(M'(\|\nabla u_k\|^{2})\right)^2\|\Delta u_{k}\|^{4}\|\nabla u_k\|^2+\frac{1}{2}\|\nabla u'_k\|^2 \notag \\
&&\qquad\leq \frac{1}{2}\beta^2\|\nabla u_k\|^{4\alpha+2}\|\Delta u_{k}\|^{4}+\frac{C_2}{2} \notag \\
&&\qquad\leq \frac{\beta^2C_2^{2\alpha+3}+C_{2}}{2}.\label{450.6}
\end{eqnarray*}
From the generalized H\"{o}lder inequality and Sobolev embedding theorem $H^2_0(\Omega)\hookrightarrow H^1_0(\Omega)\hookrightarrow L^{2(\rho+1)}(\Omega)$, we get
\begin{eqnarray*}
	\int_\Omega|u'_k|^{\rho}\nabla u''_k\nabla u'_kdx &\leq& \|u'_k\|^\rho_{2(\rho+1)}\|\nabla u'_k\|_{2(\rho+1)}\|\nabla u''_k\|  \notag \\
	&\leq& C_{s}^{\rho+1}\|\nabla u'_k\|^\rho\|\Delta u'_k\|\|\nabla u''_k\|.
\end{eqnarray*}
From Young inequality and (\ref{1023}), we deduce
\begin{equation*} 
\int_\Omega|u'_k|^{\rho}\nabla u''_k\nabla u'_kdx\leq\eta\|\nabla u''_k\|^2+\frac{C_s^{2(\rho+1)}C_{2}^{\rho}}{4\eta}\|\Delta u'_{k}\|^2.
\end{equation*}
Similary, Young inequality and assumption $(A4)$ leads to
\begin{equation}\label{452}
-\mu_2\int_\Omega \nabla u'_k\nabla z_k(x,1,t)g'(z_k(x,1,t))\leq \eta \|\nabla z_k(.,1,t)\|^2+\frac{(\mu_2c_1)^2C_{2}}{4\eta }.
\end{equation}
Taking  account to (\ref{450.6})-(\ref{452}) into (\ref{w15}) yields
\begin{eqnarray}
\nonumber&&\frac{1}{2}\frac{d}{dt}\Big[\left(1-\int_0^th(s)ds\right)\|\nabla \Delta u_k\|^2
+\|\Delta u'_k\|^2+M(\|\nabla u_k\|^2)\|\Delta u_{k}\|^{2}\\
\nonumber&&\qquad+h(\square\nabla\Delta u_{k})(t)+2\int_\Omega|u'_k|^{\rho}|\nabla u'_k|^2dx+\tau\|\nabla z_k(t)\|^2_{L^2(\Omega\times(0,1))}\Big]\\
&&\label{w16}\quad+(\frac{1}{2}-\eta)\|\nabla z_k(.,1,t)\|^2+\mu_1\int_\Omega |\nabla u'_k |^2g'(u_k)dx\\
\nonumber&&\leq\eta\|\nabla u''_k\|^2+\frac{C_s^{2(\rho+1)}C_{2}^{\rho}}{4\eta}\|\Delta u'_{k}\|^2-\frac{1}{2}h(t)\|\nabla\Delta u_{k}\|^{2}\\
\nonumber&&\qquad+\frac{1}{2}(h'\square\nabla\Delta u_k)(t)+C_{2}(\eta).
\end{eqnarray}
Multiplying (\ref{3.001}) by $g''_{ik}$ and summing with respect to  $i$, it holds that
\begin{equation*} 
\begin{split}
&\int_\Omega|u'_k|^{\rho}|u''_k|^2dx+\|\nabla u''_k\|^2\\
&=-\int_\Omega u''_k\Delta^2 u_kdx+\int_0^th(t-s)\int_{\Omega}\Delta u_k(x,s)\Delta u''_k(x,t)dxds\\
&-\int_{\Omega}M(\|\nabla u_k\|^2)\nabla u_k \nabla u''_kdx-\mu_1\int_\Omega u''_k g(u'_k)dx-\mu_2\int_\Omega u''_k g(z^k(x,1,t))dx.
\end{split}
\end{equation*}
Differentiating (\ref{3.002}) with respect to $t$, multiplying by $h'_{jk}$ and summing with respect $i$, it follows
\begin{equation}\label{w20}
\frac{\tau}{2} \frac{d}{dt}\|z'_k\|^2
+\frac{1}{2}\frac{d}{d\varrho}\|z'_k\|^2=0,
\end{equation}
Integrating (\ref{w20}) with respect to $\varrho\in(0,1)$ and summing with (3.28), we obtain that
{\small \begin{equation}\label{w21}
\begin{split}
&\int_\Omega|u'_k|^{\rho}|u''_k|^2dx+\|\nabla u''_k\|^2+\frac{\tau}{2}\frac{d}{dt}\|z'_k\|^2_{L^2(\Omega\times(0,1))}+\frac{1}{2}\|z'_k(1,t)\|^2\\
&=-\int_\Omega u''_k\Delta^2 u_kdx+\int_0^th(t-s)\int_{\Omega}\Delta u_k(x,s)\Delta u''_k(x,t)dxds+\frac{1}{2}\|u''_k\|^2\\
&-\int_{\Omega}M(\|\nabla u_k\|^2)\nabla u_k \nabla u''_kdx-\mu_1\int_\Omega u''_k g(u'_k)dx-\mu_2\int_\Omega u''_k g(z^k(x,1,t))dx.
\end{split}
\end{equation}}
In what follows, we will estimate the right hand side in (\ref{w21}). Using Green formula and Young's inequality, we get
\begin{eqnarray}
\int_\Omega u''_k(x,t)\Delta^2 u_k(x,t)dx&=&-\int_\Omega \nabla u''_k(x,t) \nabla\Delta u_k(x,t)dx \notag \\
&\leq& \eta\|\nabla u''_k(t)\|^2+\frac{1}{4\eta}\|\nabla \Delta u_k(t)\|^2. \label{w2A2}
\end{eqnarray}
Applying Cauchy-Schwarz inequality and Young's inequality, we obtain from assumption $(A2)$ and (\ref{1023}) that
\begin{eqnarray}
-\int_{\Omega}M(\|\nabla u_k\|^2)\nabla u_k \nabla u''_kdx &\leq&\delta\|\nabla u_{k}\|^{2\gamma+1}\|\nabla u''_k\| \notag \\
&\leq&\eta\|\nabla u''_k\|^2+\frac{\delta^2C_2^{2\gamma+1}}{4\eta}. \label{wxcv}
\end{eqnarray}
Similarly, we have
{\small \begin{equation*}
\begin{split}
&\int_0^th(t-s)\int_{\Omega}\Delta u_k(x,s)\Delta u''_k(x,t)dxds\\
&\quad\leq\eta \|\nabla u''_k\|^2+\frac{1}{4\eta}\int_{\Omega}\Big(\int_0^th(t-s)|\nabla\Delta u_k(s)|ds\Big)^2dx\\
&\quad\leq \eta \|\nabla u''_k\|^2+\frac{1}{4\eta}\int_{\Omega}\Big( \int_0^th(t-s)\left(|\nabla\Delta u_k(s)-\nabla\Delta u_k(t)|+|\nabla\Delta u_k(t)|\right)ds\Big)^2dx\\
&\quad=\eta \|\nabla u''_k\|^2+\frac{1}{4\eta}I.
\end{split}
\end{equation*}}
Applying H\"{o}lder's inequality and Young's, we get
{\small \begin{equation*}
\begin{split}
|I|&\leq\left(\int_0^t h(s)ds\right)^2\|\nabla\Delta u_k(t)\|^2+\int_{\Omega}\left(\int_{0}^{t}h(t-s)|\nabla\Delta u_k(s)-\nabla\Delta u_k(t)|ds\right)^2dx \\
&+2\int_{\Omega}\left(\int_0^th(t-s)|\nabla\Delta u_k(s)-\nabla\Delta u_k(t)|ds\right)\left(\int_0^th(t-s)|\nabla\Delta_ku(t)|ds\right)dx \\
&\leq 2\left(\int_0^t h(s)ds\right)^2\|\nabla\Delta_ku(t)\|^2+2\int_{\Omega}\left(\int_{0}^{t}h(t-s)|\nabla\Delta u_k(s)-\nabla\Delta u_k(t)|ds\right)^2dx.\\
&\leq 2(1-\beta_1)^2\|\nabla\Delta_ku(t)\|^2+2(1-\beta_1)(h\square\nabla\Delta u_k)(t),
\end{split}
\end{equation*}}
then, we obtain estimation
\begin{equation}\label{w22}
\begin{split}
&\int_0^th(t-s)\int_{\Omega}\Delta u_k(x,s)\Delta u''_k(x,t)dxds\\
&\leq\eta \|\nabla u''_k\|^2+\frac{(1-\beta_1)^2}{2\eta}\|\nabla\Delta_ku(t)\|^2+\frac{(1-\beta_1)}{2\eta}(h\square\nabla\Delta u_k)(t).
\end{split}
\end{equation}
And also by Young's inequality and Sobolev embedding theorem, we obtain
\begin{equation}\label{21mlk}
-\mu_1\int_\Omega u''_kg(u'_k)dx \leq\eta C_s^{2}\|\nabla u''_k\|^2+\frac{\mu_1^2}{4\eta}\int_{\Omega } |g(u'_k)^2dx.
\end{equation}
\begin{equation}\label{ert}
-\mu_2\int_\Omega u''_kg(z_k(x,1,t))dx\leq \eta C_s^{2}\|\nabla u''_k\|^2+\frac{\mu_2^2}{4\eta}\int_{\Omega }|g(z^k(x,1,t))|^2dx.
\end{equation}
Taking into account (\ref{w2A2})-(\ref{ert}) into (\ref{w21}) yields
\begin{equation}\label{w201}
\begin{split}
&\int_\Omega|u'_k|^{\rho}|u''_k|^2dx+\left(1-3\eta-2\eta C_s^{2}-\frac{C_s^{2}}{2}\right)\|\nabla u''_k(t)\|^2\\
&\qquad+\frac{\tau}{2}\frac{d}{dt}\|z'_k(t)\|^2_{L^2(\Omega\times(0,1))}+\frac{1}{2}\|z'_k(1,t)\|^2\\
&\leq\frac{2(1-\beta_1)^2+1}{4\eta}\|\nabla \Delta u_k(t)\|^2+\frac{\mu_1^2}{4\eta}\int_{\Omega } |g(u'_k)|^2dx
\\
&\qquad+\frac{\mu_2^2}{4\eta}\int_{\Omega}|g(z_k(x,1,t))|^2dx+\frac{\delta^2C_2^{2\gamma+1}}{4\eta}+\frac{1-\beta_1}{2\eta}(h\square\nabla\Delta u_k)(t).
\end{split}
\end{equation}
Thus from (\ref{w16}), (\ref{w201}), we obtain
\begin{equation}\label{w18}
\begin{split}
&\frac{1}{2}\frac{d}{dt}\Big[\left(1-\int_0^th(s)ds\right)\|\nabla \Delta u_k\|^2
+\|\Delta u'_k\|^2+M(\|\nabla u_k\|^2)\|\Delta u_{k}\|^{2}\\
&\quad+(h\square\nabla\Delta u_{k})(t)+\tau\|z'_k\|^2_{L^2(\Omega\times(0,1))}\\
&\quad+2\int_\Omega|u'_k|^{\rho}|\nabla u'_k|^2dx+\tau\|\nabla z_k(t)\|^2_{L^2(\Omega\times(0,1))}\Big]
+(\frac{1}{2}-\eta)\|\nabla z_k(.,1,t)\|^2\\
&\quad+\mu_1\int_\Omega |\nabla u'_k |^2g'(u_k)dx+\int_\Omega|u'_k|^{\rho}|u''_k|^2dx\\
&\quad+\left(1-4\eta-2\eta C_s^{2}-\frac{C_s^{2}}{2}\right)\|\nabla u''_k\|^2+\frac{1}{2}\|z'_k(1,t)\|^2\\
&\leq\frac{2(1-\beta_1)^2+1}{4\eta}\|\nabla \Delta u_k(t)\|^2+\frac{\mu_1^2}{4\eta}\int_{\Omega } |g(u'_k)|^2dx
+\frac{\mu_2^2}{4\eta}\int_{\Omega}|g(z_k(x,1,t))|^2dx\\
&\quad+\frac{C_s^{2(\rho+1)}C_{2}^{\rho}}{4\eta}\|\Delta u'_{k}\|^2+\frac{1-\beta_1}{2\eta}(h\square\nabla\Delta u_k)(t)+\frac{1}{2}(h'\square\nabla\Delta u_k)(t)\\
&\quad-\frac{1}{2}h(t)\|\nabla\Delta u_{k}\|^{2}+\frac{\delta^2C_2^{2\gamma+1}}{4\eta}+C_{2}(\eta).
\end{split}
\end{equation}
Using (A4), (\ref{1023}) and the mean value theorem, we obtain
\begin{equation} \label{w202}
\begin{array}{lll}
\int_{Q_T}|g(u'_k)|^2dxds
&=&\int_{Q_T}|g(u'_k)-g(0)||g(u'_k)|dxds\\
&\leq& c_1\int_{Q_T}g(u'_k)u'_kdxds\leq c_1C_2,
\end{array}
\end{equation}
\begin{align} \label{w2200}
\text{and}~~\int_{Q_T}|g(z_k(x,1,s))|^2dxds\leq c_1C_2,~~~\text{where}~~Q_T=\Omega\times(0,T).
\end{align}
Integrating (\ref{w18}) over $ (0,T) $, using (\ref{w202}), (\ref{w2200}), and $(A3)$, it yields
\begin{equation}\label{cw18}
\begin{split}
&\left(1-\int_0^th(s)ds\right)\|\nabla \Delta u_k\|^2+\|\Delta u'_k\|^2+M(\|\nabla u_k\|^2)\|\Delta u_{k}\|^{2}\\
&\quad+(h\square\nabla\Delta u_{k})(t)+\tau\|z'_k\|^2_{L^2(\Omega\times(0,1))}+2\int_\Omega|u'_k|^{\rho}|\nabla u'_k|^2dx\\
&\quad+\tau\|\nabla z_k(t)\|^2_{L^2(\Omega\times(0,1))}
+(1-2\eta)\int_0^t\|\nabla z_k(.,1,s)\|^2ds\\
&\quad+\int_0^t\|z'_k(.,1,s)\|^2ds+2\int_0^t\int_\Omega|u'_k|^{\rho}|u''_k|^2dxds\\
&\quad+\left(2-8\eta-4\eta C_s^{2}-C_s^{2}\right)\int_0^t\|\nabla u''_k\|^2ds+2\mu_1\int_0^t\int_\Omega |\nabla u'_k |^2g'(u_k)dxds\\
&\leq\frac{2(1-\beta_1)^2+1}{2\eta}\int_0^t\|\nabla \Delta u_k\|^2ds+\frac{\mu_1^2c_1C_2}{2\eta}
+\frac{\mu_2^2c_1C_2}{2\eta}\\
&\quad+\frac{C_s^{2(\rho+1)}C_{2}^{\rho}}{2\eta}\int_0^t\|\Delta u'_{k}\|^2ds\\
&\quad+\frac{(\delta^2C_2^{2\gamma+1}+4C_2(\eta))T}{2\eta}+\frac{1-\beta_1}{\eta}\int_0^t(h\square\nabla\Delta u_k)(s)ds+C_3,
\end{split}
\end{equation}
where $C_3$ is positive constant depending on $\|u_0\|_{H^3\cap H_0^2(\Omega)}$, $\|u_1\|_{H_0^2}$ and $\|f_0\|_{H^1_0(\Omega\times(0,1))}$. Taking $\eta$ suitably small in (\ref{cw18}) and using Gronwall lemma, we obtained the second estimate.
\begin{equation}\label{ccw18}
\begin{split}
&\|\nabla \Delta u_k(t)\|^2+\|\Delta u'_k(t)\|^2+\|z'_k(t)\|^2_{L^2(\Omega\times(0,1))}\\
&\quad+\|\nabla z_k(t)\|^2_{L^2(\Omega\times(0,1))}+\int_0^t\|\nabla z_k(.,1,s)\|^2ds\\
&\quad+\int_0^t\|\nabla u''_k(s)\|^2ds+\int_0^t\|z'_k(.,1,s)\|^2ds\leq C_4,
\end{split}
\end{equation}
where $C_4$ is positive constant depending on $\|u_0\|_{H^3\cap H_0^2(\Omega)}$, $\|u_1\|_{H_0^1}$, $\|f_0\|_{H^1_0(\Omega\times(0,1))}$, $\rho$, $g(0)$, $m_{0}$, $\beta_{1}$, $\tau$, and $T$. Estimate (\ref{ccw18}) implies
\begin{eqnarray}
\nonumber u_k \text{ is uniformly bounded in }L^{\infty}(0,T;H^3(\Omega)\cap H_0^2(\Omega)),\\
\nonumber u'_k \text{ is uniformly bounded in }L^{\infty}(0,T;H_0^2(\Omega)),\\
\nonumber u''_k \text{ is uniformly bounded in }L^2(0,T;H_0^1(\Omega)),\\
\nonumber z_k \text{ is uniformly bounded in }L^{\infty}(0,T;L^2((0,1);H^1_0(\Omega))), \\
\label{ew6} z'_k \text{ is uniformly bounded in }L^{\infty}(0,T;L^2(\Omega\times(0,1))), \\
\nonumber z_k(.,1,.) \text{ is uniformly bounded in }L^2(0,T;H^1_0(\Omega)),\\
\nonumber z'_k(.,1,.) \text{ is uniformly bounded in }L^2(0,T;L^2(\Omega)).
\end{eqnarray}
By (\ref{ew6}) and $z_{k\varrho}=-\tau z'_k$, then
\begin{equation*} 
z_k \text{ is bounded in }L^{\infty}(0,T;H^1_0(\Omega\times(0,1))).
\end{equation*}
Applying Dunford-Pettis theorem, we infer that there exists a subsequence $(u_j)$ of $(u_k)$ and $u$ such that
\begin{gather} \label{ew10}
u_k\rightharpoonup u~~ \text{ weakly star in }L^{\infty}(0,T;H^3(\Omega)\cap H_0^2(\Omega)),\\
\label{ew11}
u'_k\rightharpoonup u'~~\text{ weakly star in }L^{\infty}(0,T;H_0^2(\Omega)),\\
\label{ew12}
u''_k\rightharpoonup u''~~\text{ weakly in }L^2(0,T;H_0^1(\Omega)).
\end{gather}
By Aubin's lemma, it follows  from (\ref{ew10})-(\ref{ew12}) that there exist a subsquence of $u_j$ still denote by $u_j$ such that
\begin{equation}\label{ew13}
u_j\rightarrow u~~\text{ strongly in }L^2(0,T;H_0^2(\Omega))
\end{equation}
\begin{equation}\label{ew14}
u'_j\rightarrow u'~~\text{ strongly in }L^2(0,T;H_0^1(\Omega))
\end{equation}
which implies $\nabla u_j\rightarrow\nabla u$, $\Delta u_j\rightarrow\Delta u$ and $u'_j\rightharpoonup u'$ almost everywhere in $\Omega\times(0,T)$. Hence,
\begin{equation}\label{ew15}
|u'_j|^{\rho}u'_j\rightarrow|u'|^{\rho}u'~~\text{ almost everywhere in }\Omega\times(0,T).
\end{equation}
\begin{equation}\label{ew16}
M(\|\nabla u_j\|^2)\Delta u_j\rightarrow M(\|\nabla u\|^2)\Delta u~~\text{ almost everywhere in }\Omega\times(0,T)
\end{equation}
On the other hand, by the Sobolev embedding theorem and the first estimate, this yields
\begin{equation}\label{ew17}
\begin{array}{lll}
\||u'_j|^{\rho}u'_j\|_{L^2(0,T;L^2(\Omega))}&=&\int_0^T\|u'_j\|^{2(\rho+1)}_{2(\rho+1)}dt\\
&\leq& c_{s}^{2(\rho+1)}\int_0^T\|\nabla u'_j\|^{2(\rho+1)}_{2}dt\\
&\leq& c_{s}^{2(\rho+1)}C_{2}^{2(\rho+1)}T.
\end{array}
\end{equation}
Thus using (A2), (\ref{1023}), (\ref{w202}), (\ref{ew15}), (\ref{ew16}) and the Lions lemma (\cite{Lio},p 12.), we derive
\begin{equation}\label{ew18}
|u'_j|^{\rho}u'_j\rightharpoonup|u'|^{\rho}u'~~\text{ weakly in } L^2(0,T;L^2(\Omega)).
\end{equation}
\begin{equation}\label{ew19}
M(\|\nabla u_j\|^2)\Delta u_j\rightharpoonup M(\|\nabla u\|^2)\Delta u~~\text{ weakly in  }L^2(0,T;L^2(\Omega))
\end{equation}

\begin{equation}\label{ew20}
g(u_j)\rightharpoonup g(u)~~\text{ weakly in }L^2(0,T;L^2(\Omega)).
\end{equation}
Similarly, by applying the Dunford-Pettis theorem, we infer that there exists a subsequence $(z_j)$ of $(z_k)$ and $z$ such that
\begin{gather} \label{ew21}
z_k\rightharpoonup z~~ \text{ weakly star in }L^{\infty}(0,T;H^1_0(\Omega\times(0,1))),\\
\label{ew22}
z'_k\rightharpoonup z'~~\text{ weakly star in }L^{\infty}(0,T;L^2(\Omega\times(0,1))),\\
\nonumber z_{k\varrho}\rightharpoonup z_{\varrho}~~\text{ weakly in }L^2(0,T;L^2(\Omega\times(0,1))),\\
\label{ew24}
z_k(.,1,.) \text{ is bounded in }L^2(0,T;H^1_0(\Omega)),\\
\nonumber z'_k(.,1,.) \text{ is bounded in }L^2(0,T;L^2(\Omega)).
\end{gather}
By Aubins lemma, it follows from (\ref{ew17}), (\ref{ew18}), (\ref{ew20}) and (\ref{ew21}) that there exists a subsequence of $z_j$ still denoted by $z_j$ and subsequence of $z_j(.,1,.)$ still denoted by $z_j(.1,.)$, such that
\begin{equation*} 
z_k\rightarrow z~\text{ strongly in }L^{2}(0,T;L^2(\Omega\times(0,1))),
\end{equation*}
\begin{equation*} 
z_k(.,1,.)\rightarrow z(.,1,.)~\text{ strongly in }L^{2}(0,T;L^2(\Omega))
\end{equation*}
Also by (\ref{w2200}), (\ref{ew24}) and Lions lemma, then
\begin{equation*} 
g(z_k(.,1,.))\rightharpoonup g(z(.,1,.))~\text{ weakly in }L^{2}(0,T;L^2(\Omega)).
\end{equation*}
Let $\mathcal{D}(0,T)$ be the space of $\mathcal{C}^{\infty}$ functions with compact support in $(0,T)$. Multiplying (\ref{3.001}), (\ref{3.002}) by $\theta\in\mathcal{D}(0,T)$ and integrating over $(0,T)$, it follows that
\begin{equation}\label{ew29}
\begin{split}
&-\frac{1}{\rho+1}\int_0^T(\|u'_{k}(t)\|^{\rho}u'_k(t),w_i)\theta'(t)dt+\int_0^T(\Delta u_k(t),\Delta w_i)\theta(t)dt\\
&+\int_0^T(\nabla u''_k(t),\nabla w_i)\theta(t)dt+\int_0^TM(\|\nabla u_k(t)\|^2)(\nabla u_k(t),\nabla w_i)\theta(t)dt\\
&+\int_0^T\int_0^th(t-s)(\nabla\Delta u_k(t),\nabla w_i)\theta(t)dsdt\\
&+\mu_{1}\int_0^T(g(u'_k(t)),w_i)\theta(t)dt+\mu_{2}\int_0^T(g(z_k(.,1,t)),w_i)\theta(t)dt=0,
\end{split}
\end{equation}
{\small \begin{equation}\label{ew30}
\text{and}~~~~~~~~~\tau\int_0^T(z'_k,\varphi_i)_{L^2(\Omega\times(0,1))}\theta(t)dt+\int_0^T(z_{k\varrho},\varphi_i)_{L^2(\Omega\times(0,1))}\theta(t)dt=0.
\end{equation}}
Noting that $\{w_i\}_{i=1}^{\infty}$ is basis of $H^3(\Omega)\cap H^2_0(\Omega)$ and $\{\varphi_i\}_{i=1}^{\infty}$ is basis of $L^2(\Omega\times(0,1))$, we can pass to the limit in (\ref{ew29}), (\ref{ew30}) and obtain
\begin{equation*} 
\begin{split}
&|u'|^{\rho}u''+\Delta^{2}u-\Delta u''-M(\|\nabla u\|^{2})\Delta u-\int_{0}^{t}h(t-s)\Delta^2 u(s)ds\\
&+\mu_{1}g(u')+\mu_{2}g(z(.,1,.))=0, ~~\text{in}~~L^2(0,T;H^{-1}(\Omega)),
\end{split}
\end{equation*}
\begin{equation*} 
\tau z'+z_{\varrho}=0, ~~\text{in}~~L^2(0,T;L^2(\Omega\times(0,1)))
\end{equation*}
for arbitrary $T>0$. From (\ref{ew10})- (\ref{ew11}), (\ref{ew21})- (\ref{ew22})  and lemma 3.3.7 in \cite{ZHeng}, we conclude
$u_j(0)\rightarrow u(0)$ weakly in $H^2_0(\Omega)$,~$u'_j(0)\rightharpoonup u'(0)$ weakly in $H^1_0(\Omega)$ and $z_j(0)\rightharpoonup z(0)$ weakly in $L^2(\Omega\times(0,1))$.
Hence by (\ref{441}), we have $u(0)=u_0$, $u'(0)=u_1$ and $z(0)=f_0$. Consequently, the global existence of weak solutions is established.\\\\
\textbf{Step 2}. Uniform decay of the energy.\\
To continue our proof, we need to introduce three new functionals
\begin{equation*}
\Phi(t)=\frac{1}{\rho+1}\int_{\Omega}|u'|^{\rho}u'udx+\int_{\Omega}\nabla u'\nabla udx.
\end{equation*}
\begin{equation}\label{df4}
\begin{split}
&\Psi(t)=-\int_{\Omega}\nabla u'\int_0^th(t-s)(\nabla u(t)-\nabla u(s))dsdx\\
&~~~~~~~~-\frac{1}{\rho+1}\int_{\Omega}|u'|^{\rho}u'\int_0^th(t-s)(u(t)-u(s))dsdx.
\end{split}
\end{equation}
\begin{equation*}
\Upsilon(t)=\int_{\Omega}\int_{0}^1e^{-2\tau\varrho}G(z(x,\varrho,t))d\varrho dx
\end{equation*}
We set
\begin{equation}\label{df6}
F(t)=NE(t)+\epsilon_{1}\Phi(t)+\Psi(t)+\epsilon_2\Upsilon(t),
\end{equation}
where $N$, $\epsilon_{1}$ and $\epsilon_2$ are suitable positive constants to be determined later.
\begin{prop}
	There exist positive numbers $k_0$ and $k_{1}$ such that
	\begin{equation}\label{df7}
	k_0E(t)\leq F(t)\leq k_1E(t).\\
	\end{equation}
\end{prop}
Proof. Using (\ref{2.13}), we get
\begin{equation}\label{df8}
|\Upsilon(t)|\leq\frac{1}{\xi}E(t).
\end{equation}
Thanks to Young inequality and the Sobolev embedding theorem, we deduce
\begin{equation}\label{df9}
\begin{split}
&|\Phi(t)|\leq\frac{1}{\rho+1}\left|\int_{\Omega}|u'|^{\rho}u'udx\right|+\left|\int_{\Omega}\nabla u'\nabla u dx\right|\\
&~~~~~~~\leq\frac{1}{\rho+2}\|u'\|^{\rho+2}_{\rho+2}+\frac{(\rho+1)^{-1}}{(\rho+2)}\|u\|^{\rho+2}_{\rho+2}+\frac{1}{2}\|\nabla u'\|^{2}+\frac{1}{2}\|\nabla u\|^2\\
&~~~~~~~\leq\frac{1}{\rho+2}\|u'\|^{\rho+2}_{\rho+2}+\left(\frac{(\rho+1)^{-1}}{(\rho+2)}C_s^{\rho+2}(2E(0)/\beta_1)^{\rho/2}\right)\|\nabla u'\|^{2}\\
&\qquad+\frac{1}{2}\|\nabla u\|^2.
\end{split}
\end{equation}
By Youngs inequality and Sobolev embedding theorem, the second term in the right hand side (\ref{df4}) can be estimated as follows
\begin{equation}\label{df10}
\begin{array}{ll}
&\left|\frac{1}{\rho+1}\int_{\Omega}|u'|^{\rho}u'\int_0^th(t-s)(u(t)-u(s))dsdx\right|\\
&\leq\frac{1}{\rho+2}\|u'\|^{\rho+2}_{\rho+2}+\frac{(\rho+1)^{-1}}{(\rho+2)}\int_{\Omega}\left(\int_0^th(t-s)|u(t)-u(s)|ds\right)^{\rho+2}dx\\
&\leq\frac{1}{\rho+2}\|u'\|^{\rho+2}_{\rho+2}+\frac{(\rho+1)^{-1}}{(\rho+2)}\left(\int_0^th(s)ds\right)^{\rho+1}\int_0^th(t-s)\\
&\qquad\qquad\times\int_{\Omega}|u(t)-u(s)|^{\rho+2}dxds\\
&\leq\frac{1}{\rho+2}\|u'\|^{\rho+2}_{\rho+2}+\frac{(\rho+1)^{-1}}{(\rho+2)}(1-\beta_1)^{\rho+1}\\
&\qquad\qquad\times C_s^{\rho+2}(8E(0)/\beta_1)^{\rho/2}(h\square\Delta u)(t).
\end{array}
\end{equation}
Thus, from (\ref{df10}) we obtain
\begin{equation}\label{df11}
\begin{split}
&\Psi(t)\leq\frac{1}{2}\|\nabla u'\|^2+\frac{(1-\beta_1)C_s^2}{2}(h\square\Delta u)(t)+\frac{1}{\rho+2}\|u'\|^{\rho+2}_{\rho+2}\\
&~~~~~~+\frac{(\rho+1)^{-1}}{(\rho+2)}(1-\beta_1)^{\rho+1}C_s^{\rho+2}(8E(0)/\beta_1)^{\rho/2}(h\square\Delta u)(t).
\end{split}
\end{equation}
From (\ref{df8}),(\ref{df9}), (\ref{df11}) and the choice of $\epsilon_{1}$, $\epsilon_2$ and $N$, (\ref{df7}) can be established.\\\\
In order to obtain the exponential decay result of $E(t)$ via (\ref{cw18}), it is sufficient to prove that of $F(t)$. To this end, we need to estimate the derivative of $F(t)$ first. Using (\ref{1.1}), we obtain
\begin{equation}\label{df12}
\begin{split}
\Phi'(t)=&\frac{1}{\rho+1}\|u'\|^{\rho+2}_{\rho+2}+\|\nabla u'\|^{2}-\|\Delta u\|^2\\
&-M(\|\nabla u\|^2)\|\nabla u\|^2+\int_0^th(t-s)(\Delta u(s),\Delta u(t))ds\\
&-\mu_1\int_{\Omega}u(x,t)g(u'(x,t))dx-\mu_2\int_{\Omega}u(x,t)g(z(x,1,t))dx.
\end{split}
\end{equation}
By use of Youngs inequality and sobolev embedding theorem, we can estimate the right hand side of (\ref{df12}) as follows:
\begin{equation}\label{df13}
\begin{array}{ll}
&\int_{\Omega}\int_0^th(t-s)\Delta u(x,t)\Delta u(x,s)dsdx\\
&\leq\int_{\Omega}\int_0^th(t-s)|\Delta u(x,t)|\\
&\qquad\times(|\Delta u(x,s)-\Delta u(x,t)|+|\Delta u(x,t)|)dsdx\\
&\leq\int_0^th(s)\|\Delta u(t)\|^2ds\\
&+\int_{\Omega}\int_0^th(t-s)|\Delta u(x,t)||\Delta u(x,s)-\Delta u(x,t)|dsdx\\
&\leq(1+\eta)\int_0^th(s)\|\Delta u(t)\|^2ds+\frac{1}{4\eta}(h\square\Delta u)(t),
\end{array}
\end{equation}
\begin{equation}\label{df14}
\begin{split}
&-\mu_1\int_{\Omega}u(x,t)g(u'(x,t))dx\leq\mu_1\eta C_s^2\|\Delta u(t)\|^{2}+\frac{\mu_1}{4\eta}\|g(u'(t))\|^2,
\end{split}
\end{equation}
\begin{equation}\label{df15}
\begin{split}
-\mu_2\int_{\Omega}ug(z(x,1,t))dx\leq\mu_2\eta C_s^2\|\Delta u(t)\|^{2}+\frac{\mu_2}{4\eta}\|g(z(.,1,t))\|^2,
\end{split}
\end{equation}
where $\eta>0$. Here and in the following we use $C_s$ to represent the Poincare constant. From (A3), (\ref{df13})-(\ref{df15}), we obtain
\begin{equation}\label{df16}
\begin{split}
&\Phi'(t)\leq\frac{1}{\rho+1}\|\nabla u'\|^{\rho+2}_{\rho+2}+\|\nabla u'\|^2\\
&\qquad-\left(1-(1-\beta_1+1)(1+\eta)-(\mu_1+\mu_2)\eta C_s^2\right)\|\Delta u\|^2\\
&\leq-M(\|\nabla u\|^{2})\|\nabla u\|^2+\frac{\mu_1}{4\eta}\|g(u'(t))\|^2+\frac{\mu_2}{4\eta}\|g(z(.,1,t))\|^2\\
&\qquad+\frac{1}{4\eta}(h\square\Delta u)(t).
\end{split}
\end{equation}

Taking the derivative of $\Psi(t)$, it follows from  that (\ref{2.02})

\begin{equation}\label{df17}
\begin{array}{lll}
\Psi'(t)&=&\int_0^th(t-s)(\Delta u(t),\Delta u(t)-\Delta u(s))dsdx\\
&&-\int_0^th'(t-s)(\nabla u'(t),\nabla u(t)-\nabla u(s))ds\\
&&+\int_{\Omega}\int_0^th(t-s)M(\|\nabla u\|^2)(\nabla u(t),\nabla u(t)-\nabla u(s))dsdx\\
&&-\int_{\Omega}\left(\int_{0}^th(t-s)\Delta u(s)ds\right)\left(\int_{0}^th(t-s)(\Delta u(t)-\Delta u(s))ds\right)dx\\
&&+\mu_1\int_{\Omega}\int_0^th(t-s)(u(t)-u(s))g(u'(t))dsdx\\
&&+\mu_2\int_{\Omega}\int_0^th(t-s)(u(t)-u(s))g(z(x,1,t))dsdx\\
&&-\frac{1}{\rho+1}\int_0^th'(t-s)(|u'(t)|^\rho u'(t),u(t)-u(s))ds-\int_0^th(s)\|\nabla u'(t)\|^2ds\\
&&-\frac{1}{\rho+1}\int_0^th(s)\|u'(t)\|^{\rho+2}ds\\
&=&I_1+I_2+I_3+I_4+I_5+I_6+I_7\\
&&-\int_0^th(s)\|\nabla u'(t)\|^2ds-\frac{1}{\rho+1}\int_0^th(s)\|u'(t)\|^{\rho+2}ds.
\end{array}
\end{equation}
In what follows we will estimate $I_1,...,I_7$ in (\ref{df17}).
\begin{equation}\label{df18}
|I_1|\leq\eta\|\Delta u(t)\|^2+\frac{1-\beta_1}{4\eta}(h\square\Delta u)(t),~~\forall \eta>0.
\end{equation}
\begin{equation}\label{df19}
|I_2|\leq\eta\|\nabla u'(t)\|^2-\frac{h(0)C_s^2}{4\eta}(h'\square\Delta u)(t),~~\forall \eta>0.
\end{equation}
\begin{equation}\label{df20}
\begin{split}
&|I_3|\leq\eta(1-\beta_1)M(\|\nabla u\|^2)\|\nabla u\|^2+\frac{\delta C_s^2}{4\eta}\|\nabla u\|^{2\gamma}(h\square\Delta u)(t)\\
&~~~~\leq\eta(1-\beta_1)M(\|\nabla u\|^2)\|\nabla u\|^2+\frac{\delta C_s^2}{4\eta}(2E(0)/m_0)^{\gamma}(h\square\Delta u)(t).
\end{split}
\end{equation}
For $I_4$ in (\ref{df17}), Applying H\"{o}lder's inequality and Young's inequality, we get
\begin{equation}\label{df21}
\begin{split}
&|I_4|\leq\eta\int_{\Omega}\left(\int_0^th(t-s)\Delta u(s)ds\right)^2dx\\
&\qquad+\frac{1}{4\eta}\int_{\Omega}\left(\int_0^th(t-s)(\Delta u(t)-\Delta u(s))ds\right)^2dx\\
&\leq2\eta\left(\int_0^th(s)ds\right)^2\|\Delta u(t)\|^2\\
&\qquad+(2\eta+\frac{1}{4\eta})\left(\int_0^th(s)ds\right)(h\square\Delta u)(t),~~\forall \eta>0.
\end{split}
\end{equation}
By (A3), we obtain from (\ref{df21}) that
\begin{equation}\label{df22}
|I_4|\leq2\eta(1-\beta_1)^2\|\Delta u(t)\|^2+(2\eta+\frac{1}{4\eta})(1-\beta_1)(h\square\Delta u)(t),~~\forall \eta>0.
\end{equation}
Similarly,
\begin{equation}\label{df23}
 \begin{array}{lll}
 |I_5|&\leq&\mu_1\eta\|g(u')\|^2+\frac{\mu_1(1-\beta_1)C_s^2}{4\eta}(h\square\Delta u)(t)\\
&\leq& c_1\mu_1\eta\int_{\Omega}u'(x,t)g(u'(x,t))dx\\
&&\qquad+\frac{\mu_1(1-\beta_1)C_s^2}{4\eta}(h\square\Delta u)(t),~~\forall \eta>0.
 \end{array}
\end{equation}
\begin{equation}\label{df24}
 \begin{array}{lll}
|I_6|&\leq&\mu_2\eta\|g(z(.,1,t))\|^2+\frac{\mu_2(1-\beta_1)C_s^2}{4\eta}(h\square\Delta u)(t)\\
&\leq& c_1\mu_2\eta\int_{\Omega}z(x,1,t)g(z(x,1,t))dx\\
&&\qquad+\frac{\mu_2(1-\beta_1)C_s^2}{4\eta}(h\square\Delta u)(t),~~\forall \eta>0.
 \end{array}
\end{equation}
\begin{equation}\label{df25}
\begin{split}
|I_7|&\leq\frac{\eta}{\rho+1}\|u'(t)\|^{2(\rho+1)}_{2(\rho+1)}-\frac{h(0)C_s^2}{4\eta(\rho+1)}(h'\square\Delta u)(t)\\
&\leq\frac{\eta C_s^{2(\rho+1)}}{\rho+1}(2E(0))^{\rho}\|\nabla u'\|^2-\frac{h(0)C_s^2}{4\eta(\rho+1)}(h'\square\Delta u)(t)\\
&\leq\frac{a_0\eta}{\rho+1}\|\nabla u'\|^2-\frac{h(0)C_s^2}{4\eta(\rho+1)}(h'\square\Delta u)(t),~~\forall \eta>0,
\end{split}
\end{equation}
where $a_0=C_s^{2(\rho+1)}(2E(0))^{\rho}$. Combining (\ref{df17})-(\ref{df20}) and (\ref{df22})-(\ref{df25}) together, we arrive at
\begin{equation}\label{df26}
\begin{array}{lll}
 &&\Psi'(t)\\
 &&\leq-\frac{\int_0^th(s)ds}{\rho+1}\|u'\|^{\rho+2}_{\rho+2}+\left(\eta+\frac{a_0\eta}{\rho+1}-\int_0^th(s)ds\right)\|\nabla u'(t)\|^2\\
&& +\eta\left(1+2(1-\beta_1)^2\right)\|\Delta u(t)\|^2+c_1\mu_1\eta\int_{\Omega}u'(x,t)g(u'(x,t))dx\\
&& +c_1\mu_2\eta\int_{\Omega}z(x,1,t)g(z(x,1,t))dx\\
&& \left[\left(2\eta+\frac{1}{2\eta}+\frac{(\mu_1+\mu_2)C_s^2}{4\eta}\right)(1-\beta_1)+\frac{\delta C_s^2}{4\eta}(2E(0)/m_0)^{\gamma}\right](h\square\Delta u)(t)\\
&& -\frac{(\rho+2)h(0)C_s^2}{4(\rho+1)\eta}(h'\square\Delta u)(t),~~\forall \eta>0
\end{array}
\end{equation}
Taking also the derivative of $\Upsilon'(t)$ it follows from (\ref{2.02}) and (A4) that

\begin{equation}\label{df27}
\begin{split}
\Upsilon'(t)&=\int_{\Omega}\int_0^1e^{-2\tau\varrho}z_{\varrho}(x,\varrho,t)g(z(x,\varrho,t))\\
&=-\frac{1}{\tau}\int_{\Omega}\int_0^1\left[\frac{\partial}{\partial\varrho}(e^{-2\tau\varrho}G(z(x,\varrho,t)))+2\tau e^{-2\tau\varrho}G(x,\varrho,t)\right]d\varrho dx\\
&=-\frac{1}{\tau}\int_{\Omega}\left[e^{-2\tau}G(z(x,\varrho,t))-G(u'(x,t))\right]dx-2\Upsilon(t)\\
&\leq-2\Upsilon(t)+\frac{\alpha_2}{\tau}\int_{\Omega}u'(x,t)g(u(x,t))dx\\
&-\frac{\alpha_1e^{-2\tau}}{\tau}\int_{\Omega}z(x,1,t)g(z(x,1,t))dx.
\end{split}
\end{equation}
Then we conclude that from (\ref{df6}),(\ref{df16}), (\ref{df26}) and (\ref{df27}) that for any $t\geq t_0>0$,
\begin{equation}\label{df28}
\begin{array}{lll}
F'(t)&=&NE'(t)+\epsilon\Phi'(t)+\Psi(t)+\epsilon_2\Upsilon'(t)\\
&\leq&-\frac{h_0-\epsilon}{\rho+1}\|u'(t)\|^{\rho+2}_{\rho+2}-\left[h_0-\epsilon-\eta\left(1+\frac{a_0}{\rho+1}\right)\right]\|\nabla u'(t)\|^2\\
&&-\left[\epsilon-\eta(1-\beta_1)\right]\widehat{M}(\|\nabla u\|^2)\\
&&-\left[\epsilon\left(1-(1-\beta_1)(1+\eta)-(\mu_1+\mu_2)\eta C_s^2\right)\right. \\
&&\left. -\eta(1+2(1-\beta_1)^2)\right]\|\Delta u\|^2
\end{array}
\end{equation}
{\small \begin{equation*}
\begin{split}
&+\left[\frac{\epsilon}{4\eta}+\left(2\eta+\frac{1}{2\eta}+\frac{(\mu_1+\mu_2)C_s^2}{4\eta}\right)(1-\beta_1)+\frac{\delta C_s^2}{4\eta}(2E(0)/m_0)^{\gamma}\right](h\square\Delta u)(t)\\
&+\left[\frac{N}{2}-\frac{(\rho+2)h(0)C_s^2}{4(\rho+1)\eta}\right](h'\square\Delta u)(t)\\
&-\left(N\theta_1-\frac{\epsilon_2\alpha_2}{\tau}-\frac{\epsilon\mu_1c_1}{4\eta}-\mu_1c_1\eta\right)\int_{\Omega}u'(x,t)g(u'(x,t))dx\\
&-\left(N\theta_2+\frac{\epsilon_2\alpha_1e^{-2\tau}}{\tau}-\frac{\epsilon c_1\mu_2}{4\eta}-c_1\mu_2\eta\right)\int_{\Omega}z(x,1,t)g(z(x,1,t))dx-2\epsilon_2\Upsilon(t),~~\forall \eta>0,
\end{split}
\end{equation*}}
where $h_0=\int_0^{t_0}h(s)ds>0$, guaranteed by (A3). At this stage, we take $\epsilon<h_0$ and $\eta$ sufficiently small that
\begin{equation*}
a_2\triangleq h_0-\epsilon_1-\eta\left(1+\frac{a_0}{\rho+1}\right)>0,~~a_3\triangleq \epsilon_1-\eta(1-\beta_1)>0,
\end{equation*}
\begin{equation*}
\text{and}~~~~a_4\triangleq \epsilon_1\left(1-(1-\beta_1)(1+\eta)-(\mu_1+\mu_2)\eta C_s^2\right)-\eta\left(1+2(1-\beta_1)^2\right)>0.
\end{equation*}
Choosing $\epsilon_2>\dfrac{\xi e^{2\tau}}{2}$ for which
\begin{equation*}
-2\epsilon_2\Upsilon(t)\leq-\xi\int_{\Omega}\int_0^1G(x,\varrho,t)dxd\varrho.
\end{equation*}
As long as $\epsilon_{1}$, $\epsilon_2$ and $\eta$ are fixed, we choose $N$ large enough that
\begin{equation*}
N\theta_1-\frac{\epsilon_2\alpha_2}{\tau}-\frac{\epsilon_1\mu_1c_1}{4\eta}-\mu_1c_1\eta>0,~~N\theta_2+\frac{\epsilon_2\alpha_1e^{-2\tau}}{\tau}-\frac{\epsilon_1 c_1\mu_2}{4\eta}-c_1\mu_2\eta>0,
\end{equation*}
and
\begin{eqnarray*}
&a_5\triangleq&\xi\left[\frac{N}{2}-\frac{(\rho+2)h(0)C_s^2}{4(\rho+1)\eta}\right]-\left[\frac{\epsilon}{4\eta}+\left(2\eta+\frac{1}{2\eta}+\frac{(\mu_1+\mu_2)C_s^2}{4\eta}\right)(1-\beta_1)\right. \\
&&\qquad\left. +\frac{\delta C_s^2}{4\eta}(2E(0)/m_0)^{\gamma}\right]>0.
\end{eqnarray*}

This applying the assumption (A3) and (\ref{df28}), we deduce
\begin{equation}\label{df29}
\begin{split}
F'(t)&\leq-a_1\|u'(t)\|^{\rho+2}_{\rho+2}-a_2\|\nabla u'(t)\|^2-a_3\widehat{M}(\|\nabla u\|^2)-a_4\|\Delta u\|^2\\
&~~-a_5(h\square\Delta u)(t)-\xi\int_{\Omega}\int_0^1G(x,\varrho,t)dxd\varrho,~~\forall t\geq t_0,
\end{split}
\end{equation}
where $a_1=\frac{h_0-\epsilon}{\rho+1}$. Then (\ref{2.13}) and (\ref{df28}) imply that there exists a positive constant $M$ such that
\begin{equation}\label{df30}
F'(t)\leq-ME(t),~~\forall t\geq t_0.
\end{equation}
Combining (\ref{df7}) and (\ref{df30}), we infer
\begin{equation}\label{df31}
F'(t)\leq-\frac{M}{k_1}F(t),~~\forall t\geq t_0.
\end{equation}
Integrating (\ref{df31}) over $(t_0,t)$, it follows that
\begin{equation}\label{df32}
F(t)\leq F(t_0)e^{-\frac{M}{k_1}},~~\forall t\geq t_0.
\end{equation}
Consequently, (\ref{2.15}) can be obtained from (\ref{df7}) and (\ref{df32}). The proof is complete.

\end{document}